\documentclass[12pt]{article}
\usepackage{amsmath, amscd, amssymb,latexsym, epsfig, color}
\input xy
\xyoption{all}

\newtheorem{theorem}{Theorem}[section]

\newtheorem{lemma}[theorem]{Lemma}

\def\proof{\smallskip\noindent {\it Proof: \ }}
\def\endproof{\hfill$\square$\medskip}
 \def\Z{\mathbb{Z}}
\def\N{\mathbb{N}}

\def\m{{\mathfrak{M}}}

\def\V{{\mathbf{V}}}
\def\W{{\mathbf{W}}}
\def\S{\mathcal{S}}

\def\M{\mathcal{M}}

\newcommand{\inc}{\iota}

\newcommand{\coker}{\mbox{\upshape coker}\,}

\newcommand{\lk}{\mbox{\upshape lk}\,}
\newcommand{\cost}{\mbox{\upshape cost}\,}

\newcommand{\field}{{\bf k}}

\title{Face rings of complexes with singularities}
\author{Isabella Novik
\thanks{Research partially supported by Alfred P.~Sloan Research
Fellowship and NSF grant DMS-0801152}\\
\small Department of Mathematics, Box 354350\\[-0.8ex]
\small University of Washington, Seattle, WA 98195-4350, USA,\\[-0.8ex]
\small \texttt{novik@math.washington.edu}
\and Ed Swartz
\thanks{Research partially supported by NSF grant DMS-0900912}\\
\small Department of Mathematics, \\[-0.8ex]
\small Cornell University, Ithaca NY, 14853-4201, USA, \\[-0.8ex]
\small \texttt{ebs22@cornell.edu } }

\begin{document}
\maketitle

\begin{abstract}
It is shown that the face ring of a pure simplicial complex modulo 
$m$ generic linear forms is a ring with finite local cohomology
if and only if the link of every face of dimension $m$ or more is nonsingular. 
\end{abstract}



\section{Introduction}
In the 70's Reisner, building on unpublished work of Hochster, and 
Stanley revolutionized the study of face enumeration of simplicial complexes 
through their use of  the face ring, also called the Stanley-Reisner ring.  
Reisner proved that the face ring of a complex is Cohen--Macaulay if and only 
if the link of every face, including the empty face, is nonsingular
 \cite{Reisner}.  Here, nonsingular means that all reduced cohomology groups, 
except possibly in the maximum dimension, vanish.  Stanley used this to 
completely characterize $f$-vectors of such complexes \cite{St77}. 

A natural question which follows these results is,  
``What happens if singularities are allowed?''  The weakest relaxation 
possible is to permit nontrivial  cohomology in the lower dimensions of the 
whole complex (the link of the empty face).  Schenzel proved that for pure 
complexes the face ring is Buchsbaum if and only if this is the 
case \cite{Sch}.  The primary tool used in the proof of both Reisner's and 
Schenzel's theorems is the local cohomology of the face ring with respect to 
the irrelevant ideal.  In the Cohen--Macaulay case the local cohomology modules
 vanish below the top dimension, while in the Buchsbaum case these modules 
are finite-dimensional.  Rings with this property, that is, those  whose 
local cohomology modules below their Krull dimension are finite-dimensional, 
are called  generalized Cohen--Macaulay rings or
rings with finite local cohomology.
The goal of this note is to extend these ideas to arbitrary 
singularities.  Our main theorem says that if the dimension of the singular 
set is $m-1$, then the  the face ring modulo $m$ generic linear forms is a 
ring with finite local cohomology.   The precise statement and all definitions 
are in the next section.  This is followed by a simple example in 
Section \ref{isolated singularities} and the proof of the main theorem 
in Section \ref{proof}.

\section{Preliminaries}

For all undefined terminology we refer our readers to \cite{BrHe, St96}.
Throughout $\Delta$ is a pure $(d-1)$-dimensional simplicial complex with 
vertex set $[n]=\{1, \dots, n\}.$  If $F \in \Delta$ is a face, then the 
{\it link} of $F$ is
$$\lk F = \{G \in \Delta: F \cap G = \emptyset, F \cup G \in \Delta\}.$$
In particular, $\lk \emptyset = \Delta.$  We say that the face $F$ is 
{\it nonsingular} if $\tilde{H}^i(\lk F; \field) = 0$ for all $i < d-1-|F|.$  
Otherwise $F$ is a {\it singular} face.  The {\it singularity dimension} of 
$\Delta$ is the maximum dimension of a singular face.  If there are no 
singular faces, then $\Delta$ is {\it Cohen--Macaulay} and we  (arbitrarily) 
declare the singularity dimension of the complex to be $-\infty.$

For a field $\field$, which we will always assume is infinite (of arbitrary 
characteristic), the {\it face ring} of $\Delta$ (also known as the 
{\it Stanley-Reisner} ring) is 
$$\field[\Delta] = \field [x_1, \dots, x_n]/I_\Delta,$$
where $I_\Delta$ is the ideal generated by
$$ ( x_{i_1}\cdots x_{i_k}: \{i_1, \dots, i_k\} \notin \Delta).$$

Let $\m = (x_1, \dots, x_n)$ be the irrelevant ideal of 
$S = \field[x_1, \dots, x_n].$  For any $S$-module $M$ we use 
$H^i_\m(M)$ to denote the local cohomology modules of $M$ with respect to 
$\m$.  If the Krull dimension of $M$ is $d,$ then we say that 
$M$ is a {\it module with finite local cohomology} (or a 
{\it generalized Cohen--Macaulay} module) if for all $i<d, H^i_\m(M)$ 
is finite-dimensional as a $\field$-vector space. 
Modules with finite local cohomology were originally 
introduced in \cite{CST}, \cite{StVo} and \cite{Trung}.  
Connections between face rings and modules with finite local cohomology have been 
studied in \cite{GotoYuk} and \cite{Yuk}.  The main goal of this 
note is to prove the following.

\begin{theorem} \label{main}
Let $\Delta$ be a pure $(d-1)$-dimensional  complex.  Then the singularity 
dimension of $\Delta$ is less than $m$ if and only if for all sets of 
$m$ generic linear forms $\{\theta_1, \dots, \theta_m\}$ the quotient  
$\field[\Delta]/(\theta_1, \dots, \theta_m)$ is a 
ring  with finite local cohomology.
\end{theorem}

The proof will rely on results of Gr\"abe \cite{Grabe} which we now explain.
Denote by $|\Delta|$ the geometric realization
of $\Delta$. For a face $\tau\in \Delta$, 
let $\cost \tau :=\{\sigma \in \Delta \, : \sigma\not\supset \tau\}$
be the contrastar of $\tau$,
let $H^i(\Delta, \cost \tau)$ be the simplicial $i$-th
cohomology of the pair (with coefficients in $\field$), and for
 $\tau \subset \sigma \in \Delta$, let $\inc^*$ be the map
$H^i(\Delta, \cost \sigma) \to H^i(\Delta, \cost \tau)$ induced
by inclusion $\inc: \cost\tau \to \cost\sigma$. 
Finally, for a vector $U=(u_1, \ldots, u_n)\in \Z^n$, let 
$s(U):=\{l : u_l\neq 0\}\subseteq [n]$ be the support of $U$, let
$|U|=\sum_{l=1}^n u_l$, let $\{e_l\}_{l=1}^n$ be the standard basis 
for $\Z^n$, and let $\N$ denote the set of nonnegative integers. 

We consider the $\Z^n$-grading of $\field[x_1,\ldots, x_n]$ obtained by 
declaring $x_l$ to be of degree $e_l$. This grading refines the usual 
$\Z$-grading and induces 
a $\Z^n$-grading of $\field[\Delta]$ and its local cohomology modules.
 Thus, $H^i_\m(\field[\Delta])_{j}=
\bigoplus H^i_\m(\field[\Delta])_U$ where the sum is over all 
$U=(u_1,\ldots, u_n)\in\Z^n$ with
$|U|=j$, and multiplication by $x_l$ is a linear map from 
$H^i_\m(\field[\Delta])_U$ to $H^i_\m(\field[\Delta])_{U+e_l}$ for all 
$U\in \Z^n$.

\begin{theorem} {\rm {\bf [Gr\"abe]}}   \label{Grabe}
The following is an isomorphism of $\Z^n$-graded  $\field[\Delta]$-modules
\begin{equation}  \label{Hochster}
H^{i}_\m(\field[\Delta])\cong
\bigoplus_{\mbox{\tiny${\begin{array}{cc} -U\in \N^n\\
                             s(U)\in\Delta
           \end{array}}$}} \M^{i}_U, \quad \mbox{where} \quad
\M^{i}_U=H^{i-1}(\Delta, \cost s(U)),
\end{equation}
and the $\field[\Delta]$-structure on the $U$-th component of 
the right-hand side is given by 
$$
\cdot x_l = \left\{ \begin{array}{lll}
\text{$0$-map}, & \mbox{ if } l\notin s(U)\\
\text{identity map}, & \mbox{ if } l\in s(U) \mbox{ and } l\in s(U+e_l)\\
\inc^*: H^{i-1}(\Delta, \cost s(U)) \to H^{i-1}(\Delta, \cost s(U+e_l)), &
\mbox{ otherwise}. 
\end{array}
\right.
$$
\end{theorem}
We note that the isomorphism of (\ref{Hochster}) 
on the level of vector spaces
(rather than $\field[\Delta]$-modules) is due to Hochster, see 
\cite[Section II.4]{St96}, and that 
$$H^{i}(\field[\Delta])_0\cong \M^i_{(0,\ldots,0)} = 
H^{i-1}(\Delta, \emptyset)=\tilde{H}^{i-1}(\Delta;\field).$$

\section{Isolated singularities} \label{isolated singularities}

Before proceeding to the proof of the main theorem we consider a special case.
We say that $\Delta$ has {\it isolated singularities} if the singularity 
dimension of $\Delta$ is zero and there is at least one singular vertex.  
For the rest of the section we assume that $\Delta$ has isolated 
singularities.  

To begin with, we compute $H^i_\m(\field[\Delta])$ for $i < d.$  
Since $\Delta$ has isolated singularities Theorem \ref{Grabe} says that
$$
H^i_\m(\field[\Delta])_U = 
\begin{cases} \tilde{H}^{i-1}(\Delta; \field), & s(U)=\emptyset \\
H^{i-1}(\Delta, \cost \{j\}), & s(U) = \{j\}, \, U \in -\N^n\\
0, & \mbox{ otherwise.}
\end{cases}$$

Let $\theta =\sum^n_{i=1} a_i x_i$ be a linear form in $S$ with 
$a_i \neq 0$ for all $i.$  In order to compute 
$H^i_\m(\field[\Delta]/(\theta)$ for $i< d-1$ 
(the Krull dimension of $\field[\Delta]/(\theta)$ is $d-1$) we use the 
following short exact sequence
$$ \quad 0 \to \field[\Delta] \stackrel{\cdot \theta}{\to} 
\field[\Delta] \to \field[\Delta]/(\theta) \to 0.$$
Since multiplication by $\theta, \cdot \theta,$ is injective for any face 
ring, the above is in fact a short exact sequence.  The corresponding 
long exact sequence in local cohomology is
$$
 \dots \to H^i_\m(\field[\Delta]) \stackrel{(\cdot \theta)^\ast}{\to} 
H^i_\m(\field[\Delta]) \to 
H^i_\m(\field[\Delta]/(\theta) \stackrel{\delta}{\to} 
H^{i+1}_\m(\field[\Delta]) \stackrel{(\cdot \theta)^\ast}{\to} 
H^{i+1}_\m(\field[\Delta]) \to \dots \, .
$$
Here, $\delta$ is the connecting homomorphism which decreases the 
$\Z$-grading by one, and $\cdot \theta^\ast$ is the map induced by 
multiplication, and hence is just the module action of multiplication 
by $\theta$ on $H^i_\m(\field[\Delta]).$  
Let $f^{i}$ be the map in simplicial cohomology
$$
f^{i}:\left(\bigoplus^n_{t=1} H^{i}(\Delta, \cost \{t\})\right) 
 \to  H^{i}(\Delta, \emptyset), \,\,
f^{i}= \sum_{t=1}^n a_t \cdot \inc^\ast\left[H^{i}(\Delta, \cost \{t\})
\to  H^{i}(\Delta, \emptyset)\right].
$$
By Gr\"abe's description of the $S$-module structure of the local cohomology 
modules (using the $\Z$-grading) and the above long exact sequence,

$$H^i_\m(\field[\Delta]/(\theta))_j \cong
 \begin{cases}
 0, & j<0 \\
 \coker f^{i-1} \bigoplus \ker f^{i}, & j=0 \\
 H^i(\Delta, \emptyset), & j=1.
 \end{cases}
 $$
 In particular, $H^i_\m(\field[\Delta]/(\theta))$ is finite-dimensional.

\section{Proof} \label{proof}
We now proceed to the proof of the main theorem.
As before, let $\Delta$ be a $(d-1)$-dimensional 
simplicial complex on $[n]$. Let $\theta_1,\ldots, \theta_d$ be $d$ 
{\it generic} linear forms in $S$ with $\theta_p=\sum_{t=1}^n a_{t,p}x_t$.
In particular, we assume that every square submatrix of the 
$n\times d$-matrix $A=(a_{t,p})$ is non-singular  
and the $\theta$'s satisfy the prime avoidance argument 
in the proof of Lemma \ref{isomorphism} below. 
Each  $\theta_p$ acts on $\field[\Delta]$ by multiplication,
$\cdot\theta_p: \field[\Delta] \to \field[\Delta]$. This map, in turn, 
induces the map 
$(\cdot \theta_p)^{\ast}=\cdot \theta_p \ :H^l_{\m}(\field[\Delta])
\to H^l_{\m}(\field[\Delta])$ that increases the $\Z$-grading by one.
The key objects in the proof are the kernels of these maps and their 
intersections:
$$
\ker^l_{m, i} := \bigcap_{p=1}^m \left(\ker (\cdot \theta_p)^{\ast} \, : 
\, H^l_\m(\field[\Delta])_{-(i+1)}
\to  H^l_\m(\field[\Delta])_{-i}\right) \quad \mbox{and} \quad 
\ker^l_{m}:=\bigoplus_{i\in\Z} \ker^l_{m, i}.
$$
Thus
$\ker^l_m$ is a graded submodule of 
$H^l_\m(\field[\Delta])$ and $\ker^l_{m,i}$ 
is simply  $(\ker^l_m)_{-(i+1)}$.

What are the dimensions of these kernels? If $m=0$, 
then $\ker^l_{0, i} = H^l_\m(\field[\Delta])_{-(i+1)}$, 
and Theorem \ref{Grabe} implies that for $i\geq 0$,
\begin{eqnarray}  \nonumber
\dim_{\field} \ker^l_{0, i} & = & \sum_{F\in\Delta}
|\{U\in\N^n \, : \, s(U)=F, \, |U|=i+1\}| \cdot 
\dim_{\field} H^{l-1}(\Delta, \cost F) \\
\label{Hochs}
& = &  
\sum_{F\in\Delta}
\binom{i}{|F|-1} \cdot \dim_{\field} H^{l-1}(\Delta, \cost F).
\end{eqnarray}
For a general $m$, we prove the following (where
we set $\binom{a}{b}=0$ if $b<0$).

\begin{lemma}    \label{ineq}
For every $0\leq m \leq d$, $l\leq d$, and $i\geq m$, 
$$
\dim_{\field} \ker^l_{m,i} \leq 
\sum_{F\in\Delta}
\binom{i-m}{|F|-m-1}\cdot \dim_{\field}  H^{l-1}(\Delta, \cost F).
$$
\end{lemma}

In fact, using Lemma \ref{ineq}
(but deferring its proof to the end of the section), 
we can say even more:

\begin{lemma}  \label{equality}
For every $0\leq m \leq d$, $l\leq d$, and $i\geq m$, 
$$
\dim_\field \ker^l_{m,i} = 
\sum_{F\in\Delta}
\binom{i-m}{|F|-m-1} \cdot \dim_{\field} {H}^{l-1}(\Delta, \cost F). 
$$
Moreover, the map 
$\oplus_{i\geq m+1} \ker^l_{m,i} \stackrel{\cdot \theta_{m+1}}
{\longrightarrow} \oplus_{i\geq m+1}\ker^l_{m,i-1}$, 
is a surjection. 
\end{lemma}
\proof We prove the statement on the dimension of 
$\ker^l_{m,i}$ by induction on $m$. 
For $m=0$ (and any $l, i \geq 0$), this is eq.~(\ref{Hochs}). 
For larger $m$, we notice that
the restriction of 
$(\cdot \theta_{m+1})^{\ast}$ to 
$\ker^l_{m,i}$ is a linear map from $\ker^l_{m,i}$
to $\ker^l_{m,i-1}$, whose kernel is $\ker^l_{m+1,i}$. 
Thus for $i\geq m+1$,
\begin{eqnarray*}
\dim_\field \ker^l_{m+1,i}& \geq & 
\dim_\field \ker^l_{m,i}-\dim \ker^l_{m,i-1}\\
 &=&
\sum_{F\in\Delta}
\left[\binom{i-m}{|F|-m-1}-\binom{i-1-m}{|F|-m-1}\right] \cdot
 \dim_\field H^{l-1}(\Delta, \cost F) \\
 &=&\sum_{F\in\Delta}
\binom{i-(m+1)}{|F|-(m+1)-1}\cdot \dim_\field {H}^{l-1}(\Delta, \cost F).
\end{eqnarray*}
The second step in the above computation 
is by the inductive hypothesis. Comparing
the resulting inequality to that of Lemma \ref{ineq} 
shows that this inequality is in 
fact equality, and hence that the map 
$(\cdot\theta_{m+1})^\ast \, : \,   \ker^l_{m,i} \to \ker^l_{m,i-1}$
is surjective for $i\geq m+1$. 
\endproof

Lemma \ref{equality} allows us to get a handle on 
$H^l_{\m}(\field[\Delta]/(\theta_1, \ldots, \theta_m))_{-i}$ at least for
$l,i>0$:

\begin{lemma} \label{isomorphism}
For $0\leq m \leq d$ and $0<l\leq d-m$, there is a
graded isomorphism of modules
$$
\bigoplus_{i\geq 1} H^l_\m(\field[\Delta]/(\theta_1,\ldots, \theta_m))_{-i} 
\cong \bigoplus_{i\geq1} \ker^{l+m}_{m,i+m-1} \; \mbox{with} \;
H^l_\m(\field[\Delta]/(\theta_1,\ldots, \theta_m))_{-i} \cong
 \ker^{l+m}_{m,i+m-1}.  
$$ 
Thus, for $l,i>0$,
$$\dim_\field H^l_\m(\field[\Delta]/(\theta_1,\ldots, \theta_m))_{-i}= 
\sum_{F\in\Delta}
 \binom{i-1}{|F|-m-1}\cdot \dim_\field {H}^{l+m-1}(\Delta, \cost F).$$ 
\end{lemma}

\proof The proof is by induction on $m$, with the $m=0$ case being evident.
For larger $m$, we want to mimic the proof given in 
Section \ref{isolated singularities}.
One obstacle to this approach is that the map 
$\cdot\theta_{m+1}: \field[\Delta]/(\theta_1,\ldots, \theta_m)\to
\field[\Delta]/(\theta_1,\ldots, \theta_m)=:\M[m]$ 
might not be injective anymore. 
However, a ``prime avoidance'' argument 
together with the genericity assumption on $\theta_{m+1}$ implies that
$\cdot\theta_{m+1}:\M[m]/H^0_\m(\M[m]) \to \M[m]/H^0_\m(\M[m])$ is injective,
and hence one has the corresponding long exact sequence in local cohomology.   See, for instance, \cite[Chapter 3]{Eis} for details on ``prime avoidance" arguments.  
Since $H^0_\m(\M[m])$ has Krull dimension 0, modding $H^0_\m(\M[m])$ out 
does not affect $H^l_\m$ for $l>0$, 
so that the part of this long exact
sequence corresponding to $l,i>0$ can be rewritten as
\begin{eqnarray*}
\bigoplus_{i\geq 1} H^l_\m(\M[m])_{-(i+1)} 
 &\stackrel{(\cdot \theta_{m+1})^\ast}{\longrightarrow}&
\bigoplus_{i\geq 1} H^l_\m(\M[m])_{-i}  \longrightarrow
\bigoplus_{i\geq 1} H^l_\m(\M[m+1])_{-i} \\  &\stackrel{\delta}{\longrightarrow}& 
\bigoplus_{i\geq 1} H^{l+1}_\m(\M[m])_{-(i+1)} 
\stackrel{(\cdot \theta_{m+1})^\ast}{\longrightarrow} 
\bigoplus_{i\geq 1} H^{l+1}_\m(\M[m])_{-i}.
\end{eqnarray*}
By the inductive hypothesis combined with Lemma \ref{equality}, the leftmost
map in this sequence is surjective. Hence the module
$\oplus_{i\geq 1} H^l_\m(\M[m+1])_{-i}$ is
isomorphic to the kernel of the rightmost map.
Applying the inductive hypothesis to the last two entries of 
the sequence then implies the following isomorphism of modules:
$$\bigoplus_{i\geq 1} H^l_\m(\M[m+1])_{-i} \cong  
\bigoplus_{i\geq 1} \ker ((\cdot \theta_{m+1})^\ast \, : \, 
\ker^{l+m+1}_{m, i+m} \to \ker^{l+m+1}_{m, i+m-1}) = 
\bigoplus_{i\geq 1} \ker^{l+m+1}_{m+1, i+m}.$$
\endproof

Theorem \ref{main} now follows easily from Lemma \ref{isomorphism}:

\smallskip\noindent {\it Proof of Theorem \ref{main}: \ }
Since $\field[\Delta]$ is a finitely-generated algebra,
$H^0_\m(\field[\Delta]/(\theta_1,\ldots, \theta_m))$ 
has Krull dimension zero, and hence is a finite-dimensional 
vector space for any simplicial complex $\Delta$. 
So we only need to care about $H^l_\m$ for $l>0$.
As $\binom{i-1}{|F|-m-1}>0$ for all $i\gg 0$ and $|F|>m$, 
Lemma \ref{isomorphism} implies that 
$\field[\Delta]/(\theta_1,\ldots, \theta_m)$ is a ring with
finite local cohomology if and only if for all faces
$F\in\Delta$ of size larger than $m$ and all  $l+m<d$, 
the cohomology 
${H}^{l+m-1}(\Delta, \cost F)$ vanishes. Given that 
 ${H}^{l+m-1}(\Delta, \cost F)$ is isomorphic to $
\tilde{H}^{l+m-1-|F|}(\lk F; \field)$ 
(see e.g. \cite[Lemma 1.3]{Grabe}),
this happens if and only if each such $F$ is non-singular.
\endproof

To finish the proof of Theorem \ref{main}, it only remains to verify 
Lemma \ref{ineq}. For its proof we use the following notation. 
Fix $m, l >0$, and $i\geq m$. 
For $r\in \{i, i+1\}$, let 
$\V_r:=\{U\in\N^n \, : \, |U|=r, \, s(U)\in\Delta\}$, 
and for $F\in\Delta$, let $\V_{r,F}:= \{u\in\V_r \, : \, s(U)=F\}$.
If $F=\{f_1<\ldots<f_j\}\in\Delta$ where $j > m$, then set
$
\W_{r, F} := \{U=(u_1,\ldots, u_n)\in \V_{r,F} \, : \, 
u_{f_s}=1 \mbox{ for } 1\leq s\leq m \}.
$
Observe that $\W_{i+1, F}$ is a subset of $\V_{i+1, F}$ of cardinality 
$\binom{i-m}{|F|-m-1}$.

For $G\subseteq F\in\Delta$ let 
$\Phi_{F,G}$ denote the map $\inc^\ast \, : \, H^{l-1}(\Delta, \cost F) 
\to H^{l-1}(\Delta, \cost G)$.
Thus, $\Phi_{F,G}$ is the identity map if $F=G$.
Using Theorem \ref{Grabe}, we identify $H^l_\m(\field[\Delta])_{-r}$
with $\bigoplus_{U\in\V_{r}} H^{l-1}(\Delta, \cost s(U))$, 
and for $z\in H^l_\m(\field[\Delta])_{-r}$ 
we write $z=(z_U)_{U\in\V_{r}}$, 
where $z_U\in H^{l-1}(\Delta, \cost s(U))$ is the $(-U)$-th component
of $z$. 
Since $\theta_p=\sum_{t=1}^n a_{t,p}x_t$, Theorem~\ref{Grabe} yields
that for such $z$, $r=i+1$, and $T\in \V_i$,
\begin{equation} \label{mult}
(\theta_p z)_T = \sum_{\{t \, : \, T+e_t\in\V_{i+1}\}}
    a_{t,p} \cdot \Phi_{s(T+e_t), s(T)}(z_{T+e_t}).
\end{equation}

\smallskip\noindent {\it Proof of Lemma \ref{ineq}: \ } 
Since $|\W_{i+1, F}|=\binom{i-m}{|F|-m-1}$ for all $F\in\Delta$,
to prove that
$\dim_{\field} \ker^l_{m,i} \leq \sum_{F\in\Delta}
\binom{i-m}{|F|-m-1}\cdot \dim_{\field}  H^{l-1}(\Delta, \cost F),$
it is enough to verify that for $z,z'\in \ker^l_{m,i}$, 
$$  z_U=z'_U \ \mbox{ for all }  F\in\Delta \mbox{ and all } 
U\in\W_{i+1, F} \quad 
\Longrightarrow \quad z=z',
$$ 
or equivalently (since $\ker^l_{m,i}$ is a $\field$-space) that
for $z\in \ker^l_{m,i}$,
\begin{equation}  \label{initial}
z_U=0 \ \mbox{ for all }  F\in\Delta \mbox{ and all } 
U\in\W_{i+1, F} \quad 
\Longrightarrow \quad z=0.
\end{equation}
To prove this, fix such a $z$. From eq.~(\ref{mult}) 
and the definition of $\ker^l_{m,i}$, it follows that 
\begin{equation}  \label{zero}
\sum_{\{t \, : \, T+e_t\in\V_{i+1}\}}
    a_{t,p} \cdot \Phi_{s(T+e_t), s(T)}(z_{T+e_t}) =0 
\quad \forall \ 1\leq p \leq m \mbox{ and } \forall \  T\in \V_i.
\end{equation}
For a given $T\in\V_i$, we refer to the $m$ conditions imposed on $z$
by eq.~(\ref{zero}) as ``the system defined by $T$", and denote this 
system by $\S_T$.

Define a partial order, $\succ$, 
on $\V_{i}$ as follows: $T'\succ T$ if 
either $|s(T')|>|s(T)|$, or $s(T')=s(T)$ and the last non-zero
entry of $T'-T$ is positive. To finish the proof, we
verify by a descending (w.r.t~$\succ$)
induction on $T\in\V_i$, that $z_{T+e_t}=0$ for all $t$ with
$T+e_t\in\V_{i+1}$. For $T\in\V_i$, there are two possible 
cases (we assume that $s(T)=F=\{f_1< \ldots <f_j\}$):

\smallskip\noindent{\bf Case 1:}  $T\in\W_{i,F}$. 
Then for each $t\notin F$, either $F':=F\cup\{t\}\notin\Delta$,
in which case $T+e_t\notin\V_{i+1}$, or $T+e_t\in\W_{i+1, F'}$,
in which case $z_{T+e_t}=0$ by eq.~(\ref{initial}). Similarly,
if $t=f_r$ for some $r>m$, then $T+e_t\in \W_{i+1, F}$, and 
$z_{T+e_t}=0$ by (\ref{initial}). Finally, 
for any $t\in F$, $s(T+e_t)=s(T)=F$, and so $\Phi_{s(T+e_t), s(T)}$
is the identity map. Thus the system $\S_T$ 
reduces to $m$ linear equations in $m$ variables:
\begin{equation} \label{system}
\sum_{r=1}^m a_{f_r, p}\cdot z_{T+e_{f_r}} =0  
\quad \forall  \  1\leq p \leq m.
\end{equation}
Since the matrix $(a_{f_r,p})_{1\leq r,p \leq m}$ is non-singular, 
$(z_{T+e_t}=0 \mbox{ for all }t)$ is the only solution of $\S_T$.

\medskip\noindent{\bf Case 2:} $T\notin \W_{i,F}$, and so $T_{f_s}\geq 2$
for some $s\leq m$. Then for any $t\notin F$, $T':=T-e_{f_s}+e_t \succ T$,
as $T'$ has a larger support than $T$, and $T'+e_{f_s}=T+e_t$. Hence 
$z_{T+e_t}=z_{T'+e_{f_s}}=0$ by the inductive hypothesis on $T'$. 
Similarly, if $t\in F-\{f_1, \ldots, f_m\}$, then $t>f_s$, 
and so $T'':= T-e_{f_s}+e_t \succ T$. As $T''+e_{f_s}=T+e_t$, 
the inductive hypothesis on $T''$ imply that $z_{T+e_t}=0$. Thus,
as in Case 1, $\S_T$ reduces to system (\ref{system}),
whose only solution is trivial. 
\endproof


\end{document}